\input amssym.def
\input amssym
\magnification=1200
\parindent0pt
\hsize=16 true cm
\baselineskip=13  pt plus .2pt
$ $

\def\G{(\Gamma,{\cal G})}

\def\D{{\Bbb D}}

\def\H{{\cal H}}

\centerline {\bf A note on the Nielsen realization problem for connected
sums of $S^2 \times S^1$}

\bigskip

\centerline {Bruno P. Zimmermann}

\centerline {Universit\`a degli Studi di Trieste}

\centerline {Dipartimento di Matematica e Geoscienze}

\centerline {34127 Trieste, Italy}

\bigskip  \bigskip

{\bf Abstract.}  {\sl   We consider finite group-actions on 3-manifolds
$\H_g$ obtained as the connected sum of $g$ copies of $S^2
\times S^1$, with free fundamental group $F_g$ of rank $g$.  We prove that,
for $g > 1$,  a  finite group of diffeomorphisms  of
$\H_g$ inducing a trivial action on homology is cyclic and
embeds into an $S^1$-action on $\H_g$. As a consequence, no
nontrivial element of the twist subgroup of the mapping class group of
$\H_g$  (generated by Dehn twists
along embedded 2-spheres) can be realized by a periodic diffeomorphism of
$\H_g$ (in the sense of the Nielsen realization problem). We also discuss
when a finite subgroup of the outer automorphism group 
${\rm Out}(F_g)$ of the fundamental group of $\H_g$ can be realized by a
group of diffeomorphisms of $\H_g$.}

\medskip


{\sl Key words and phrases:}  3-manifold; connected sums of $S^2
\times S^1$; finite group action; mapping class group; outer automorphism
group of the fudamental group; Nielsen realization problem

\bigskip

{\bf 1. Introduction}

\medskip

All finite group-actions in the present paper will  be faithful, smooth and
orientation-preserving, all manifolds orientable.  We
are interested in finite group-actions on connected sums
$\H_g = \sharp_g (S^2 \times S^1)$ of $g$ copies of $S^2 \times S^1$; we will
call $\H_g$ a  {\it closed handle of genus $g$} in the
following. The fundamental group of $\H_g$ is the free group
$F_g$ of rank $g$. Considering induced actions on the fundamental
group and on the first homology $H_1(\H_g) \cong \Bbb Z^g$, there are   
canonical maps

$${\rm Diff}(\H_g) \to {\rm Out}(F_g) \to {\rm GL}(g, \Bbb Z)$$

where ${\rm Diff}(\H_g)$ denotes the orientation-preserving diffeomorphism
group of  $\H_g$ and
${\rm Out}(F_g) = {\rm Aut}(F_g)/{\rm Inn}(F_g)$ the outer automorphism
group of its fundamental group.

\bigskip

{\bf Theorem 1.}  {\sl  Let $G$ be a finite group acting on a closed handle
$\H_g$ of genus $g > 1$ such that the induced action on the first homology of
$\H_g$ is trivial. Then $G$ is  
cyclic and a subgroup of an $S^1$-action on $\H_g$; in particular,
all elements of $G$ are isotopic to the identity.}

\medskip

For a description and
classification of circle-actions on 3-manifolds and closed handles, see 
[14].

\medskip

Denoting by  ${\rm Mod}(\H_g)$ the {\it mapping class group} of isotopy
classes of orientation-preser-ving diffeomorphisms of $\H_g$, there are
induced maps

$${\rm Mod}(\H_g) \to {\rm Out}(F_g) \to {\rm GL}(g, \Bbb Z).$$

Let ${\rm Twist}(\H_g)$ denote the subgroup of ${\rm Mod}(\H_g)$ generated
by all Dehn twists along embedded 2-spheres in $\H_g$  (i.e., by   
cutting along a 2-sphere and regluing after twisting by one full turn
around an axis; since such a twist represents a generator of $\pi_1(SO(3))
\cong \Bbb Z_2$, its square is isotopic to the identity). 
By classical results of Laudenbach [6],[7] there is a short exact sequence

$$ 1  \to  {\rm Twist}(\H_g)  \hookrightarrow {\rm Mod}(\H_g) \to  {\rm
Out}(F_g) \to 1;$$

moreover ${\rm Twist}(\H_g) \cong  (\Bbb Z_2)^g$ is generated by the sphere
twists around the core spheres $S^2 \times *$ of the $g$ different $S^2
\times S^1$ summands of $\H_g$ (twists around seperating
2-spheres instead are isotopic to the identity).
It is proved in [1] that ${\rm Mod}(\H_g)$
is isomorphic to a semidirect product  ${\rm Twist}(\H_g) \rtimes
{\rm Out}(F_g)$.  Theorem 1 has the following consequence (in the sense
of the {\it Nielsen realization problem}).

\bigskip

{\bf Corollary 1.}  {\sl  No nontrivial element of the 
twist group  ${\rm Twist}(\H_g)$ can be realized (represented) by a
periodic diffeomorphism of $\H_g$.}

\bigskip

For $g > 1$ this follows from Theorem 1 but the methods apply
also to  the case $g = 1$ of $\H_1 = S^2 \times S^1$, using the
fact that $S^2 \times S^1$ is a geometric
3-manifold belonging to the ($S^2 \times \Bbb R$)-geometry (one of
Thurston's eight 3-dimensional geometries, see [15]), and that finite
group-actions on $S^2 \times S^1$ are geometric ([10, Theorem 8.4]).

\medskip

For a  solution of the Nielsen realization problem for aspherical and Haken
3-manifolds, see [21] (here finite groups of mapping classes can always be realized,
except for a purely algebraic obstruction in the case of Seifert fiber spaces where,
however, a finite inflation of the group can always be realized).

\medskip

By [6], homotopic diffeomorphisms of $\H_g$ are isotopic 
but this does not
remain true for arbitrary connected sums of 3-manifolds. By [4], 
twists around separating 2-spheres in  a 3-manifold may or may not be homotopic to the
identity, moreover by [3] there are sphere-twists which are homotopic but
not isotopic to the identity (see also the discussion in the
introduction of [1]).  As an example, considering a connected sum
$M = M_1 \sharp M_2$ of two closed hyperbolic 3-manifolds $M_1$ and $M_2$, 
the sphere-twist around the connecting 2-sphere is not
homotopic to the identity; also, it cannot be realized by a periodic map
(e.g., if $M_1$ or $M_2$ does not admit a nontrivial periodic map
then also the connected sum $M = M_1 \sharp M_2$ has no periodic maps).

\medskip

There arises naturally the question of which finite subgroups
of  ${\rm Out}(F_g)$  can be realized by a finite group of diffeomorphisms
of $\H_g$.  Finite groups $G$ of diffeomorphisms of
$\H_g$ which act faithfully on the fundamental group (i.e., inject into
${\rm Out}(F_g)$) are considered in [17] where, for $g
\ge 15$, the quadratic upper bound
$|G| \le 24g(g-1)$  for their orders is obtained. Since  ${\rm Out}(F_g)$
has finite subgroups of  larger orders, these subgroups cannot be
realized by finite groups of diffeomorphisms (by [16] the
maximal order of a finite subgroup of 
${\rm Out}(F_g)$ is $2^gg!$, for $g > 2$). A precise result is as follows
(we refer to [17, section 2] for definitions and the proof).

\bigskip

{\bf Theorem 2.}  {\sl  Let $G$ be a finite subgroup of 
${\rm Out}(F_g)$  and $1 \to F_g \to E \to  G \to 1$ the corresponding group
extension associated to $G$. Then $G$ can be realized by an isomorphic group
of diffeomorphisms of $\H_g$  if and only if $E$ is
isomorphic to the fundamental group $\pi_1\G$ of a finite graph of finite
groups  $\G$ in normal form associated to a closed handle-orbifold
(in particular, the vertex groups $\G$ have to be isomorphic to finite subgroups
of SO(4) and the edge groups to finite subgroups of SO(3)). }

\bigskip

We note that, for a finite group $G$ acting on a closed handle $\H_g$, the
quotient $\H_g/G$ has the structure of a closed handle-orbifold (see [17]).
Analogous results on finite group-actions on 3-dimensional handlebodies
are obtained in [8] and [12] (and in [9] for finite group-actions on
handlebodies in arbitrary dimensions).

\medskip

The case $g = 2$ is special. By well-known
results, 
$${\rm Out}(F_2) \cong {\rm Aut}(\Bbb Z^2)  \cong {\rm GL}(2, \Bbb Z) 
\cong  \D_6 *_{\D_2}\D_4,$$ 
so up to conjugation the maximal finite subgroups of ${\rm Out}(F_2)$ are the dihedral
groups
$\D_6$ and $\D_4$ of orders 12 and 8, and both can be realized by diffeomorphisms of the
torus with one boundary component (hence, if the realizations of the amalgamated subgroups
$\D_2$ coincide, one obtains a realization of the whole group ${\rm Out}(F_2) \cong  \D_6
*_{\D_2}\D_4$).  Considering the product with a closed interval, one obtains realizations
on the handlebody $V_2$ of genus 2 and alsos on its double $\H_2$ along the boundary.

\medskip

Concerning the case
$g = 3$,  by [18] there are exactly five  maximal finite subgroups of ${\rm Out}(F_3)$ up
to conjugation; by an easy application of Theorem 2, all of these maximal finite
subgroups can be realized by  diffeomorphisms of the closed handle $\H_3$ of genus 3 (but
not of a handlebody $V_3$ of genus 3).

\bigskip

{\bf 2. Proof of Theorem 1}

\medskip

Let $G$ be a finite group acting faithfully and orientation-preservingly on
a  closed handle
$\H_g = \sharp_g (S^2 \times S^1)$ of genus $g$.  By the equivariant sphere
theorem (see [10] for an approach by minimal surface  techniques, [2] and
[5] for topological-combinatorial proofs), there exists an embedded,
homotopically nontrivial 2-sphere  $S^2$ in $\H_g$ such that
$x(S^2) = S^2$ or
$x(S^2) \cap S^2 = \emptyset$ for all $x \in G$. We cut $\H_g$ along the system
of disjoint 2-spheres $G(S^2)$, by removing the interiors of  
$G$-equivariant regular neighbourhoods $S^2 \times [-1,1]$ of these
2-spheres, and call each of these regular neighbourhoods $S^2
\times [-1,1]$ a 1-handle. The result is a collection of 3-manifolds with
2-sphere boundaries, with an induced action of $G$. We close each of the
2-sphere boundaries by a 3-ball and extend the action of $G$ by taking the
cone over the  center of each of these 3-balls, so $G$ permutes these
3-balls and their centers. The result is a finite collection of closed
handles of lower genus on which $G$ acts (cf. [17]). Applying inductively the procedure
of cutting along 2-spheres, we finally end up with a finite collection of
3-spheres or 0-handles (closed handles of
genus 0). Note that the construction gives a finite graph $\Gamma$ on which
$G$ acts whose vertices correspond to the 0-handles and whose edges to the
1-handles. Note that $\Gamma$ has no {\it free edges}, i.e. edges with one
vertex of valence 1.

\medskip

On each 3-sphere (0-handle) there are finitely many points which are the
centers of the attached 3-balls (their boundaries are the 2-spheres along
which the 1-handles are attached).  For each of these 3-spheres, let 
$G_v$ denote its stabilizer in $G$ (by the  geometrization of finite group-actions
on 3-manifolds, one may assume that the action of a stabilizer
$G_v$  on the corresponding 3-sphere is orthogonal but this is not needed for
the following).
Denoting by $G_e$ the stabilizer in $G$ of a 1-handles $S^2
\times [-1,1]$, we can assume that each stabilizer $G_e$ preserves the
product structure of $S^2 \times [-1,1]$ of the corresponding 1-handle (by
choosing small equivariant regular neighbourhoods of the 2-spheres). If some
element of a stabilizer
$G_e$ acts as a reflection on [-1,1], we split the 1-handle into two
1-handles by introducing a new 0-handle obtained from a small regular
neighbourhood
$S^2 \times [-\epsilon,\epsilon]$ of 
$S^2 \times \{0\}$ by closing up with two 3-balls. Hence we can assume that each
stabilizer $G_e$ of a 1-handle $S^2 \times [-1,1]$ does not interchange its
two boundary 2-spheres; that is, $G$ acts {\it without inversions} on the
graph $\Gamma$.

\medskip

Suppose now that $g > 1$ and that the induced action of $G$ on the
first homology of $\H_g$ and hence also of  $\Gamma$ is
trivial. As before,
$G$ acts without inversions on $\Gamma$ and $\Gamma$ has no free edges.  We
will prove in next Proposition that under these hypotheses the action
of $G$ on $\Gamma$ is trivial, that is each element of $G$ acts as the
identity  on $\Gamma$.  Hence $G$ fixes each vertex and each edge of
$\Gamma$. 

\medskip

Since $G$ fixes each 1-handle $S^2 \times [-1,1]$, it maps each 2-sphere 
$S^2 \times \{0\}$ to itself. By construction, $G$ does not interchange the two sides ot 
such a 2-sphere and acts faithfully on it (otherwise some element of $G$ would act
trivially on  an invariant regular neighbourhood of such a 2-sphere and then act
trivially also on all of $\H_g$ (well-known in particular for smooth actions)). It follows
that
$G$ is isomorphic to a finite subgroup of the orthogonal group SO(3), i.e. cyclic $\Bbb
Z_n$, dihedral $\Bbb D_{2n}$, tetrahedral
$\Bbb A_4$, octahedral $\Bbb S_4$ or dedecahedral $\Bbb A_5$. 
It is easy to see that an orientation-preserving action of
$\Bbb D_{2n}$, $\Bbb A_4$, $\Bbb S_4$ or $\Bbb A_5$ on
$S^3$  has  at most two global fixed points around which a 1-handle can be
attached; but then the graph $\Gamma$ would be a segment or a circle, that is
$g \le 1$. Since
$g > 1$, $G$ is a cyclic group which acts by rotations around an axis $S^1$
in each 0-handle $S^3$. By the positive solution of the Smith-conjecture [13],
each of these axes is a trivial knot in $S^3$, and hence the action of
the cyclic group $G$ embeds into an
$S^1$-action on each 0-handle. Since these $S^1$-actions on the 0-handles
extend to the connecting 1-handles $S^2 \times  [-1,1]$, the cyclic
$G$-action on $\H_g$ embeds into an $S^1$-action.

\medskip

To complete the proof of Theorem 1, it remains to prove the
following Proposition (which may be considered as an analogue of Theorem  1
for finite graphs).

\bigskip

{\bf Proposition.}  {\sl Let $G$ be a finite group acting faithfully on
a finite connected graph $\Gamma$ without free edges and of genus $g > 1$ (or cycle
rank, or rank of its free fundamental group). Then also
the induced action of $G$ on the first homology $H_1(\Gamma) \cong \Bbb Z^g$
of $\Gamma$ is faithful.}

\bigskip

{\it Proof.}  By subdividing edges, we can
assume that $G$ acts without inversion of edges on $\Gamma$.  Suppose  that
an element  $x \in G$  acts trivially on the first homology of $\Gamma$.
Then its Lefschetz number is $1-g$ which, by the Hopf trace formula, is
equal to the Euler characteristic of the fixed point set of $x$ which is a
subgraph  $\Gamma'$ of $\Gamma$ (since $G$ acts without
inversions of edges). The graph $\Gamma$ of genus $g$ has Euler
characteristic $1-g$; passing from $\Gamma'$ to $\Gamma$ by adding
successively the missing edges, the Euler characteristic remains unchanged (when
adding a free edge) or decreases. Since $\Gamma$ has no free edges, this
implies $\Gamma' = \Gamma$, and hence $x$ acts trivally on
$\Gamma$. This completes the proof of the Proposition.

\bigskip

By [19, proof of Satz 3.1],
each finite subgroup of ${\rm Out}(F_g)$ can be realized by an action of the
group on a finite graph $\Gamma$ without free edges (this is a version of the
{\it Nielsen realization problem for finite graphs} which several years later was 
"rediscovered" by various authors); the Proposition implies then the 
following well-known result.

\bigskip

{\bf Corollary 2.}  {\sl  The canonical projection 
${\rm Out}(F_g) \to  {\rm GL}(g, \Bbb Z)$ is injective on finite subgroups
of ${\rm Out}(F_g)$.}

\bigskip

We note that not all finite subgroups of ${\rm GL}(g, \Bbb Z)$ are induced
in this way by finite subgroups of ${\rm Out}(F_g)$; in fact, for $g = 2, 4,
6, 7, 8, 9$ and 10 there are finite subgroups of ${\rm GL}(g, \Bbb Z)$ of
orders larger than $2^gg!$ (which, by [16], is the maximal order of a finite
subgroup of ${\rm Out}(F_g)$)). On the other hand, there are also
small cyclic subgroups of ${\rm GL}(g, \Bbb Z)$ which cannot be realized in
this way, see the discussion in [20, section 5].

\bigskip  \bigskip
\vfill \eject

\centerline {\bf References}

\bigskip

\item {[1]}  T. Brendle, N. Broaddus, A. Putman, {\it  The mapping class
group of connected sums of $S^2 \times S^1$,} arXiv 2012:01529

\smallskip

\item {[2]}  M.J. Dunwoody,  {\it  An equivariant sphere theorem,} Bull.
London Math. Soc. 17  (1985), 437-448  

\smallskip

\item {[3]}  J.L. Friedman, D.M. Witt,  {\it  Homotopy is not istopy for
homeomorphisms of 3-manifolds,}  Topology 25 (1986), 35-44

\smallskip

\item {[4]}  H. Hendriks,  {\it  Applications de la th\'eorie d'obstruction
en dimension 3,} Bull. Soc. Math. France M\'em.  53  (1977), 81-196

\smallskip

\item {[5]}  W. Jaco, J.H. Rubinstein,  {\it  PL equivariant surgery and
invariant decompositions of 3-manifolds,}  Adv. in Math. 73  (1989),
149-191

\smallskip

\item {[6]} F. Laudenbach, Sur le 2-sph\`eres d'une vari\'et\'e de diemsion
3,  Ann. Math. 97 (1973), 57-81

\smallskip

\item {[7]}  F. Laudenbach,  Topologie de la dimension trois: homotopie et
isotopie,  Soci\'et\'e Math\'e-matique de France, Paris 1974

\smallskip

\item {[8]} D. McCullough, A. Miller, B. Zimmermann,  {\it Group actions
on handlebodies,}  Proc. London Math. Soc.  59   (1989), 373-415

\smallskip

\item {[9]}  M. Mecchia, B. Zimmermann,  {\it On finite groups of
isometries of handlebodies in arbitrary dimensions and finite extensions of
Schottky groups,}  Fund. Math. 230 (2015),  237-249

\smallskip

\item {[10]} W.H. Meeks, P. Scott,  {\it Finite group actions on
3-manifolds,}  Invent. math. 86  (1986) , 287-346

\smallskip

\item {[11]} W.H. Meeks, L. Simon, S.T. Yau,  {\it Embedded minimal
surfaces, exotic spheres, and manifolds with positive Ricci curvature,} 
Ann. of Math. 116  (1982) , 621-659

\smallskip

\item {[12]} A. Miller, B. Zimmermann,  {\it  Large groups of symmetries of
handlebodies,}  Proc. Amer. Math. Soc. 106  (1989),  829-838

\smallskip

\item {[13]} J.W. Morgan,  H. Bass,  {\it The Smith Conjecture,}  
Academic Press 1984

\medskip

\item {[14]} F. Raymond,  {\it Classification of actions of the circle on
3-manifolds,} Trans. Amer. Math. Soc. 131  (1968), 51-78

\smallskip

\item {[15]}  P. Scott,  {\it The geometries of 3-manifolds,} 
Bull. London Math. Soc. 15  (1983),  401-487

\smallskip

\item {[16]} S. Wang, B. Zimmermann,  {\it  The maximum order finite groups
of outer automorphisms of free groups,}  Math. Z. 216  (1994), 83-87

\smallskip

\item {[17]} B. Zimmermann,  {\it On finite groups acting on a connected sum
of 3-manifolds $S^2 \times S^1$,} Fund. Math. 226  (2014), 131-142

\smallskip

\item {[18]} B. Zimmermann,  {\it Finite groups of outer automorphism groups
of free groups,}  Glasgow Math. J. 38  (1996),  275-282

\smallskip

\item {[19]} B. Zimmermann,   {\it \"Uber Hom\"oomorphismen n-dimensionaler
Henkelk\"orper und endliche Erweiterungen von Schottky-Gruppen,}   Comm.
Math. Helv. 56   (1981), 474-486

\smallskip

\item {[20]} B. Zimmermann,  {\it  Lifting finite groups of outer
automorphisms of free groups, surface groups and their abelianizations,}
Rend. Istit. Mat. Univ. Trieste 37  (2005), 273-282  (electronic
version  under  http://rendiconti.dmi.units.it)

\smallskip

\item {[21]} B. Zimmermann,   {\it Das Nielsensche Realisierungsproblem f\"ur
hinreichend grosse 

3-Mannigfaltigkeiten,}   Math. Z.  180   (1982), 349-359

\bye